\documentclass[a4paper, 10pt]{article}
\usepackage{amsmath,amssymb}
\usepackage[all]{xy}
\usepackage{abstract}
\usepackage{authblk}

\usepackage{amsthm}
\theoremstyle{plain}
\newtheorem{theorem}{Theorem}[section]

\newtheorem{corollary}[theorem]{Corollary}
\newtheorem{lemma}[theorem]{Lemma}
\newtheorem{definition}[theorem]{Definition}

\theoremstyle{definition}
\newtheorem{example}[theorem]{Example}
\newtheorem{remark}[theorem]{Remark}
\usepackage{blindtext}
\usepackage{sectsty}
\sectionfont{\centering}
\subsectionfont{\centering}
\numberwithin{equation}{section}

\begin{document}

\title{\bfseries{On the volume functions and the cohomology rings of special weight varieties of type $A$}}
\author{Tatsuru TAKAKURA and Yuichiro YAMAZAKI}
\date{}

\maketitle

\begin{abstract}
In this paper, we consider the cohomology rings of some multiple weight varieties of type $A$, that is, symplectic torus quotients for a direct product of several coadjoint orbits of the special unitary group. Under some specific assumptions, we prove the symplectic volumes of multiple weight varieties are equal to the volumes of flow polytopes.  Using differential equations satisfied by the volume functions of flow polytopes, we give an explicit presentation of the cohomology ring of the multiple weight variety of special type.
\end{abstract}

\section{Introduction}
Let $G$ be a compact and connected Lie group with Lie algebra $\mathfrak{g}$, and $T$ a maximal torus of $G$ with Lie algebra $\mathfrak{t}$.  Let $\mathfrak{g^{*}}$ and $\mathfrak{t}^{*}$ be the dual vector spaces of $\mathfrak{g}$ and $\mathfrak{t}$, respectively. Using an invariant inner product on $\mathfrak{g}$, we identify $\mathfrak{g^{*}}$ and $\mathfrak{t^{*}}$ with $\mathfrak{g}$ and $\mathfrak{t}$, respectively. Under this identification, we  interpret $\mathfrak{t^{*}}$ as a subspace of $\mathfrak{g}^{*}$. 

Let ${\mathcal O}_{\lambda}$ be the coadjoint orbit of $G$ through $\lambda \in \mathfrak{t}^{*} \subset \mathfrak{g^{*}}$. It is well known that $\mathcal{O}_{\lambda}$ has the $G$-invariant symplectic structure called the Kostant-Kirillov-Souriau symplectic form. For $\lambda_{1}, \lambda_{2}, \ldots, \lambda_{n}, \mu \in \mathfrak{t^{*}}$, we consider the symplectic quotient
\begin{eqnarray*}
{\mathcal M}_{T}
&:=&({\mathcal O}_{\lambda_{1}} \times \cdots \times {\mathcal O}_{\lambda_{n}})/\hspace{-2pt}/_{\mu}T\\
&=& \left\{(x_{1},\ldots,x_{n})\in {\mathcal O}_{\lambda_{1}} \times \cdots \times {\mathcal O}_{\lambda_{n}}\  \middle| \ \sum_{i=1}^{n}\Phi_{\lambda_{i}}(x_{i})=\mu \right\}/T
\end{eqnarray*}
called a multiple weight variety, where  $\Phi_{\lambda_{i}}:\mathcal{O}_{\lambda_{i}} \to \mathfrak{t}^{*}$ is the moment map for the $T$ action on $\mathcal{O}_{\lambda_{i}}$, and $T$ acts diagonally on the product of coadjoint orbits.

When $n=1$, many results have been known for ${\mathcal M}_{T}$. For example, in the case $G=U(n)$, some formulas for the volume of weight varieties $\mathcal{M}_{T}$ are given in [7], and explicit presentations of the cohomology ring of $\mathcal{M}_{T}$ are given in [6].

When $n=2$ and $G=SU(3)$, a volume formula of $\mathcal{M}_{T}$ and some examples are given in [15]. 
 
Our first result is the following theorem that plays an important role in this paper.

\begin{theorem}[{see Theorem 2.1}]
Let $(M,\omega)$ be a symplectic manifold endowed with a free Hamiltonian $T$-action, and let $\Phi: M\to \mathfrak{t^{*}}$ be a proper moment map. Suppose that $\omega$ and $\Phi$ can be written as $\omega=\sum\limits_{
i=1}^{n}p_{i}\omega_{i}$ and $\Phi=\sum\limits_{i=1}^{n}p_{i}\Phi_{i}$, where $p_{i} \in \mathbb{R}$, $\omega_{i}$ is a closed form, and each $\Phi_{i}$ satisfies the condition $\iota_{\xi_{Y}}\omega_{i}=d\langle \Phi_{i}, Y \rangle$ for all $Y \in \mathfrak{t}$. We fix a regular value $\mu_{0} \in \mathfrak{t}^{*}$ of $\Phi$. Let $\omega_{\mu_{0}}$ be the symplectic form on $M_{T}=\Phi^{-1}(\mu_{0})/T$. 
Then there exist an open neighborhood $U$ of $\mu_{0}$ and a diffeomorphism $\phi_{\mu}: \Phi^{-1}(\mu_{0})/T \to \Phi^{-1}(\mu)/T$ for all $\mu \in U$. Moreover there exist the cohomology classes $z_{i}$ and $v_{j} \in H^{2}(M_{T};\mathbb{R})$ such that 
\[
\phi_{\mu}^{*}([\omega_{\mu}])=\sum_{i=1}^{n}p_{i}z_{i} + \sum_{j=1}^{m}x_{j}v_{j}
\]
for all $\mu=\sum\limits_{j=1}^{m}x_{j}X_{j}\in U$, where $x_{j}\in  \mathbb{R}$ and $X_
{1},\ldots, X_{m}$ is a basis of $\mathfrak{t}^{*}$ .
\end{theorem}

Immediately, we have the following corollary.
\begin{corollary}[{see Corollary 2.2}]
Under the same assumptions of Theorem 1.1, suppose that $H^{*}_{T}(M;\mathbb{R})$ is generated by the equivariant cohomology classes $[\omega_{i}-\Phi_{i}] \in H^{2}_{T}(M;\mathbb{R})$ $(1\leq i \leq n)$ as a $H^{*}(BT;\mathbb{R})$-module. Then the cohomology classes $z_{i}$ and $v_{j}$ in Theorem 2.1 generate the cohomology ring $H^{*}(M_{T};\mathbb{R})$ multiplicatively.
\end{corollary}

We consider the case $G=SU(l+1)$ of type $A_{l}$. Let $\alpha_{1},\ldots,\alpha_{l}$ be the simple roots and $W$ be the Weyl group. We denote by $\mathfrak{t}^{*}_{+}$ and $\mathfrak{t}^{*}_{++}$ the Weyl chamber and its interior respectively. 

When $\lambda_{1},\ldots, \lambda_{n}\in \mathfrak{t}^{*}_{++}$, then $\mathcal{O}_{\lambda_{1}}\times \cdots \times \mathcal{O}_{\lambda_{n}}$ satisfies the assumptions in Theorem 1.1 and Corollary 1.2 (see Theorem 3.1).

Let $C(A_{l}^{+})$ be the convex cone generated by $A_{l}^{+}=\left\{\alpha_{1},\ldots, \alpha_{l}  \right\}$ and we set $\rho=\frac{1}{2}\sum\limits_{i=1}^{l}\alpha_{i}$.

Let $e_{1},\ldots, e_{l+1}$ be the standard basis of $\mathbb{R}^{l+1}$. We assign a positive integer $m_{i,j}$ to each $i$ and $j$ $(1\leq i<j\leq l+1)$ and set $m=(m_{i,j})$ and $M=\sum\limits_{1\leq i<j\leq l+1} m_{i,j}$.

For $h \in C(A_{l}^{+})$, the polytope
\[
P_{l,m}(h)=
\left\{ 
(x_{i,j}^{k})\in \mathbb{R}^{M}
\ \middle| \ 
\begin{aligned}
1\leq i < j \leq l+1,\ \  1 \leq k \leq m_{i,j}, \\
 x_{i,j}^{k} \geq 0,\ \   \sum_{i<j}\sum_{k}x_{i,j}^{k}(e_{i}-e_{j})=h
\end{aligned}
\right\}
\]
is called the flow polytope.
For a positive integer $n$, we denote by $P_{l,n}(h)$ the flow polytope when $m_{i,j}=n$ for all $i,j$.

For $\lambda_{1},\ldots, \lambda_{n}, \mu \in \mathfrak{t}^{*}$, we say that $\mu$ is sufficiently close to $\lambda=\lambda_{1}+\cdots+\lambda_{n}$ when the following condition holds:

If $\sigma_{1}(\lambda_{1}+\rho)+\cdots+\sigma_{n}(\lambda_{n}+\rho)-(\mu+n\rho) \in C(A_{l}^{+})$ for $\sigma_{1},\ldots, \sigma_{n}\in W$, then $\sigma_{1}=\cdots=\sigma_{n}=e$.

When $\mu$ is sufficiently close to $\lambda$, the symplectic volume ${\rm vol}(\mathcal{M}_{T})$ of the multiple weight variety $\mathcal{M}_{T}$ is equal to the volume $v_{l,n}(\lambda-\mu)$ of the flow polytope $P_{l,n}(\lambda-\mu)$. 

\begin{theorem}[see Theorem 3.3]
Suppose that $G=SU(l+1)$, $\lambda_{1}, \ldots, \lambda_{n} \in \mathfrak{t}_{+}^{*}-\{0\}$ and let $\mu \in \mathfrak{t}^{*}$ be a regular value  of the moment map sufficiently close to $\lambda=\lambda_{1}+\cdots +\lambda_{n}$. Then
\[
{\rm vol}(\mathcal{M}_{T})=v_{l,n}(\lambda-\mu).
\]
\end{theorem}

The open subset $\mathfrak{c}_{{\rm nice}}$ of $C(A_{l}^{+})$ is defined by
\[
\mathfrak{c}_{{\rm nice}}=\left\{h:=q_{1}\alpha_{1}+\cdots+q_{l}\alpha_{l}  \  \middle| \  q_{1},\ldots , q_{l} \in \mathbb{R}_{\geq 0},\  0<q_{1}<q_{2}<\cdots <q_{l} \right\},
\]
which is called the nice chamber.
In [14], it was proved that if $h\in \mathfrak{c}_{\rm nice}$, then volume functions $v_{l,m}(h)$ of flow polytopes satisfy a certain system of differential equations and conversely, the solution of the system of differential equations is unique up to constant multiple. We prove that if $h\in \mathfrak{c}_{\rm nice}$, then Ann($v_{l,m}$)$=\{{\rm differential\  operator}\  \partial \ |$ $\ \partial v_{l,m}=0 \}$ is generated by the differential operators defined in [14] (see Theorem 3.8).

As we will see in Section 2, if $\lambda_{1},\ldots, \lambda_{n}\in \mathfrak{t}^{*}_{++}$, then there exists an one-to-one correspondence between differential equations satisfied by the volume functions and relations of the cohomology rings of multiple weight varieties. 
Thus, we have the following result.

\begin{theorem}[{see Theorem 3.9}]
Suppose that $G=SU(l+1)$ and $\lambda_{1}, \ldots, \lambda_{n} \in \mathfrak{t}^{*}_{++}$. Let $\mu \in \mathfrak{t}^{*}$ be a regular value of the moment map sufficiently close to $\lambda_{1}+\cdots+\lambda_{n}$, such that $\lambda_{1}+\cdots +\lambda_{n}-\mu$ is in the nice chamber. Then there exist cohomology classes $z_{1}, \ldots, z_{l} \in H^{2}(\mathcal{M}_{T};\mathbb{R})$  such that the cohomology ring $H^{*}(\mathcal{M}_{T};\mathbb{R})$ is given by
\[
H^{*}(\mathcal{M}_{T};\mathbb{R})\cong \frac{\mathbb{R}[z_{1},z_{2}, \ldots, z_{l}]}{(z_{l}^{n}, z_{l-1}^{n}(z_{l-1}+z_{l})^{n},\ldots, z_{1}^{n}(z_{1}+z_{2})^{n}\cdots (z_{1}+\cdots +z_{l})^{n})}.
\]
\end{theorem}

This paper is organized as follows. In Section 2, we prove Theorem 1.1, Corollary 1.2 and review the one-to-one correspondence between differential equations and relations of the cohomology rings. In Section 3, we prove Theorem 1.3. Using a characterization of the volume functions of  flow polytopes and Theorem 1.3, we prove Theorem 1.4.

\section{Symplectic volume and cohomology}
\subsection{The symplectic form on a torus quotient}
In this subsection we prove a generalization of Theorem ${\rm 9.8.1}$ in [10].

Let $T$ be a torus with Lie algebra $\mathfrak{t}$ and let $\mathfrak{t}^{*}$ be the dual of $\mathfrak{t}$. Using an invariant inner product, we identify  $\mathfrak{t}$ with $\mathfrak{t}^{*}$.

Let $(M, \omega)$ be a symplectic manifold endowed with a symplectic $T$-action. An action of $T$ on $(M, \omega)$ is said to be Hamiltonian if there exists a $T$-equivariant map $\Phi:M \to \mathfrak{t}^{*}$ called a moment map that satisfies
\[
\iota_{\xi_{X}}\omega=d \langle \Phi, X \rangle
\]
where $\xi_{X}$ is the fundamental vector field associated with $X \in \mathfrak{t}$ defined by
\[
(\xi_{X})_{x} =\frac{d}{dt}({\rm exp}tX \cdot x)|_{t=0}
\]
for $x \in M$.

Suppose that $\mu$ is a regular value of $\Phi$, that $\Phi^{-1}(\mu)$ is not empty and that the $T$-action on $\Phi^{-1}(\mu)$ is free. Then $\Phi^{-1}(\mu)$ is a smooth submanifold of $M$ and $\Phi^{-1}(\mu)/T$ is a smooth symplectic manifold called a symplectic quotient. In this paper, we consider the case that $\mu \in \mathfrak{t}^{*}$ satisfies the conditions as above. 

If an action of $T$ on $M$ is free, then $\pi : M \to M/T$ is a principal $T$-bundle. In this situation, equipping this bundle with a connection, we have the horizontal subspace $\Omega(M)_{{\rm hor}}$, curvature forms $\Omega_{i}\in \Omega(M)_{{\rm hor}}$ and a $T$-equivariant map
\begin{align}
S(\mathfrak{t}^{*})\otimes \Omega(M)\to \Omega(M)_{{\rm hor}}\hspace{15pt}  X_{I}\otimes \eta \mapsto \Omega_{I}\wedge \eta_{{\rm hor}}
\end{align}
where $S(\mathfrak{t}^{*})$ is the symmetric algebra on $\mathfrak{t}^{*}$, $X=(X_{1},\ldots,X_{m})$ is a basis of $\mathfrak{t}^{*}$ and $I=(i_{1},\ldots, i_{m})$ is a multi-index. The subcomplex $\Omega(M)_{{\rm bas} }:=\pi^{*}\Omega(M/T)\subset \Omega(M)$ is called  the complex of basic forms. A form $\omega$ is basic if and only if it is $T$-invariant and horizontal. 

As explaind in [10], we have a bijection $\pi^{*}: \Omega(M/T)\to \Omega(M)_{{\rm bas}}$.
Combining the map (2.1), restricted to invariant forms, with the inverse of $\pi^{*}$, we have a map
\[
C:\Omega_{T}(M):=(S(\mathfrak{t}^{*})\otimes \Omega(M))^{T}\to \Omega(M/T)
\]
called the Cartan map ([10]).

Based on the proof of Theorem 9.8.1 in [10], we have the followiong theorem that plays an important role in Section 3.

\begin{theorem}
Let $(M,\omega)$ be a symplectic manifold endowed with a free Hamiltonian $T$-action, and let $\Phi: M\to \mathfrak{t^{*}}$ be a proper moment map. Suppose that $\omega$ and $\Phi$ can be written as $\omega=\sum\limits_{
i=1}^{n}p_{i}\omega_{i}$ and $\Phi=\sum\limits_{i=1}^{n}p_{i}\Phi_{i}$ where $p_{i} \in \mathbb{R}$, $\omega_{i}$ is a closed form, and each $\Phi_{i}$ satisfies the condition $\iota_{\xi_{Y}}\omega_{i}=d\langle \Phi_{i}, Y \rangle$ for all $Y \in \mathfrak{t}$. We fix $\mu_{0} \in \mathfrak{t}^{*}$ that satisfies the conditions as above. Let $\omega_{\mu_{0}}$ be the symplectic form on $M_{T}=\Phi^{-1}(\mu_{0})/T$. 
Then there exist an open neighborhood $U$ of $\mu_{0}$ and a diffeomorphism $\phi_{\mu}: \Phi^{-1}(\mu_{0})/T \to \Phi^{-1}(\mu)/T$ for all $\mu \in U$. Moreover there exist the cohomology classes $z_{i}$ and $v_{j} \in H^{2}(M_{T};\mathbb{R})$ such that 
\[
\phi_{\mu}^{*}([\omega_{\mu}])=\sum_{i=1}^{n}p_{i}z_{i} + \sum_{j=1}^{m}x_{j}v_{j}
\]
for all $\mu=\sum\limits_{j=1}^{m}x_{j}X_{j}\in U$  where $x_{j}\in  \mathbb{R}$ and $X_
{1},\ldots, X_{m}$ is a basis of $\mathfrak{t}^{*}$.
\end{theorem}

\begin{proof}
Since $\omega=\sum\limits_{
i=1}^{n}p_{i}\omega_{i}$ and each $\omega_{i}$ is a closed form, we have
\[
[\omega]=\sum_{i=1}^{n}p_{i}[\omega_{i}] \in H^{2}(M;\mathbb{R}).
\]
As explained in [10], since $\Phi=\sum\limits_{i=1}^{n}p_{i}\Phi_{i}$ and each $\Phi_{i}$ satisfies the condition $\iota_{\xi_{Y}}\omega_{i}=d\langle \Phi_{i}, Y \rangle$, the equivariant 2-form $\omega-\Phi$, known as the equivariant symplectic form, is closed and so
\[
[\omega-\Phi]=\sum_{i=i}^{n}p_{i}[\omega_{i}-\Phi_{i}] \in H^{2}_{T}(M;\mathbb{R}).
\]

Since $\mu_{0}$ is a regular value of $\Phi$, there exists an open neighborhood $U$ of $\mu$ such that $\Phi^{-1}(U)\to U$ is a trivial bundle, and so we have a diffeomorphism $\phi_{\mu}: \Phi^{-1}(\mu_{0})/T \to \Phi^{-1}(\mu)/T$ for all $\mu \in U$. 

For any $\mu \in U$,  we consider the following commutative diagram.

\vspace{7pt}
\[
\xymatrix{
\Phi^{-1}(\mu)  \ar[r]^{i_{\mu}}  \ar[d]_{\pi_{\mu}} & M  \ar[d]_{\pi} \ar[r]^-\Phi &\mathfrak{t}^{*}   \\
\Phi^{-1}(\mu)/T \ar[r]^{j_{\mu}} &M/T \ar[r]^-\Psi & \mathfrak{t}^{*} \ar@{=}[u]
}
\]
Let $\omega_{\mu}$ be the symplectic form on $\Phi^{-1}(\mu)/T$ which satisfies
$\pi^{*}_{\mu}\omega_{\mu}=i^{*}_{\mu}\omega$
and note that the map $\Psi$ is defined by
$\Phi=\Psi \circ \pi$.

For the equivariant symplectic form $\omega-\Phi$, we have
\[
c:=C(\omega -\Phi)=C\left(\omega-\sum\Phi_{X_{j}}X_{j}\right)=\nu-\sum\Psi_{j}\Omega_{j},
\]
where $C$ is the Cartan map and $\nu$ is the unique form on $M/T$ with the property 
\[
\pi^{*}\nu=\omega_{\rm{hor}}\in \Omega(M)_{\rm bas}.
\]
Since $i_{\mu}^{*}\omega \in \Omega(\Phi^{-1}(\mu))_{{\rm bas}}$, we have $i_{\mu}^{*}\omega=i_{\mu}^{*}\omega_{{\rm hor}}$,
and so
\[
\pi_{\mu}^{*}j_{\mu}^{*}\nu=i_{\mu}^{*}\pi^{*}\nu=i_{\mu}^{*}\omega_{{\rm hor}}=\pi_{\mu}^{*}\omega_{\mu}.
\]
Since $\pi_{\mu}^{*}$ is injective, we have $j_{\mu}^{*}\nu=\omega_{\mu}$.

Therefore we have
\begin{eqnarray*}
\phi^{*}_{\mu}([\omega_{\mu}])&=&\phi^{*}_{\mu} \circ j_{\mu}^{*}\left([c]+\sum\limits_{j=1}^{m}\Psi_{j}[\Omega_{j}]\right)\\
&=&\sum\limits_{i=1}^{n} p_{i} \phi^{*}_{\mu} \circ j_{\mu}^{*}\circ C([\omega_{i}-\Phi_{i}])+\sum\limits_{j=1}^{m} x_{j}\phi^{*}_{\mu} \circ j_{\mu}^{*}\circ C([X_{j}])\\
&=&\sum\limits_{i=1}^{n} p_{i}z_{i}+\sum\limits_{j=1}^{m} x_{j}v_{j},
\end{eqnarray*}
where $z_{i}$ and $v_{j}$ are defined by $z_{i}=\phi^{*}_{\mu} \circ j_{\mu}^{*} \circ C([\omega_{i}-\Phi_{i}])$, $v_{j}=\phi^{*}_{\mu} \circ j_{\mu}^{*}\circ C([X_{j}])$.
\qedhere
\end{proof}

The map $j_{\mu_{0}}^{*}\circ C : \Omega_{T}^{*}(M)\to \Omega^{*}(M_{T})$ in Theorem 2.1
induces the surjective map called the Kirwan map ([11])
\[
\kappa:H^{*}_{T}(M;\mathbb{R})\to H^{*}(M_{T};\mathbb{R}).
\]

\begin{corollary}
Under the same conditions as those in Theorem 2.1, suppose that $H^{*}_{T}(M;\mathbb{R})$ is generated by the equivariant cohomology classes $[\omega_{i}-\Phi_{i}] \in H^{2}_{T}(M;\mathbb{R})$ $(1\leq i \leq n)$ as a $H^{*}(BT;\mathbb{R})$-module. Then the cohomology classes $z_{i}$ and $v_{j}$ in Theorem 2.1 generate the cohomology ring $H^{*}(M_{T};\mathbb{R})$ multiplicatively.
\end{corollary}
\begin{proof}
Since the equivariant cohomology ring $H^{*}_{T}(M;\mathbb{R})$ is generated by the equivariant cohomology classes $[\omega_{i}-\Phi_{i}] \in H^{2}_{T}(M;\mathbb{R})$ $(1\leq i \leq n)$ as a $H^{*}(BT;\mathbb{R})$-module and $z_{i}=\kappa([\omega_{i}-\Phi_{i}])$, $v_{j}=\kappa([X_{j}])$, we have the result.
\end{proof}

Let $(M,\omega)$ be a compact symplectic manifold. Recall that the symplectic volume ${\rm{vol}}(M)$ is defined by
\[
{\rm{vol}}(M)=\int_{M} e^{\omega}.
\]
In the case of Theorem 2.1, the symplectic volume ${\rm vol}(M_{T})$ is a polynomial of $p_{i}\  (1\leq i \leq n)$ and  $x_{j}\   (1\leq j \leq m)$ and is written as
\begin{align}
&{\rm vol}(M_{T})=\int_{M_{T}}\frac{1}{d!}(p_{1}z_{1}+\cdots+p_{n}z_{n}+x_{1}v_{1}+\cdots+x_{m}v_{m})^{d} \notag \\
&=\int_{M_{T}}\sum_{d_{1},\ldots,d_{n+m}}\frac{1}{d_{1}! \cdots d_{n+m}!}(p_{1}z_{1})^{d_
{1}}\cdots (x_{m}v_{m})^{d_{
n+m}}\notag \\
\label{vol}
&=\sum_{d_{1},\ldots,d_{n+m}}\frac{p_{1}^{d_{1}}\cdots p_{n}^{d_{n}}x_{1}^{d_{n+1}}\cdots x_{m}^{d_{n+m}}}{d_{1}!\cdots d_{n+m}!}\int_{M_{T}}z_{1}^{d_{1}}\cdots z_{n}^{d_{n}}v_{1}^{d_{n+1}}\cdots v_{m}^{d_{n+m}}
\end{align}
where $d_{1},\ldots, d_{n+m}\in \mathbb{Z}_{\geq 0}$ and $\sum\limits_{i=1}^{n+m}d_{i}=d:=\frac{1}{2}{\rm dim}_{\mathbb{R}}M_{T}$.

Immediately, we have the following.

\begin{corollary}
Under the same conditions as those in Theorem 2.1 and Corollary 2.2, the symplectic volume ${\rm vol}(M_{T})$ is the generating function of intersection pairings
\[
\int_{M_{T}}z_{1}^{d_{1}}\cdots z_{n}^{d_{n}}v_{1}^{d_{n+1}}\cdots v_{m}^{d_{n+m}}
\]
 of $M_{T}$, where $d_{1},\ldots, d_{n+m}\in \mathbb{Z}_{\geq 0}$ and $\sum\limits_{i=1}^{n+m} d_{i}=\frac{1}{2}{\rm {dim}}_{\mathbb{R}}M_{T}$ .
\end{corollary}

\subsection{Volume functions and Poincar$\acute{{\text e}}$ duality algebras}

We review the definition of a Poincar$\acute{{\text e}}$ duality algebra based on [1].

\begin{definition}
Let $\Bbbk$ be a field. Let $\mathcal{A}^{*}=\bigoplus_{j=0}^{d}\mathcal
{A}^{2j}$ be a finite dimensional graded commutative $\Bbbk$-algebra such that
\begin{enumerate}
\item there exists an isomorphism $\int_{\mathcal{A}}: \mathcal{A}^{2d}\to \Bbbk$,
\item the pairing $\mathcal{A}^{2p} \otimes \mathcal{A}^{2d-2p}\to \Bbbk, a \otimes b \mapsto \int_{\mathcal{A}}a\cdot b$ is non-degenerate.
\end{enumerate}
Then $\mathcal{A}$ is called a Poincar$\acute{{\text e}}$ duality algebra of formal dimension $2d$.
\end{definition}

Consider the algebra of differential operators with constant coefficients $\mathcal{D}:=\mathbb{R}[\partial_{1},\ldots,\partial_{n}]$, where $\partial_{i}=\frac{\partial}{\partial t_{i}}$ for $1\leq i \leq n$. We assume that deg$\partial_{i}=2$. For any non-zero homogeneous polynomial $v \in \mathbb{R}[t_{1},\ldots, t_{n}]$ of degree $d$, we consider the following ideal in $\mathcal{D}$
\[
{\rm Ann}(v):=\{D\in \mathcal{D}\  |\ Dv=0\}.
\]
It is known that the quotient $\mathcal{D}/{\rm Ann}(v)$ is a Poincar$\acute{{\text e}}$ duality algebra of formal dimension $2d$.

Based on the proof of Theorem 9.8.2 in [10], we have the following lemma.

\begin{lemma}
Let $(M,\omega)$ be a compact symplectic manifold. Suppose that $w_{1},\ldots,$ $w_{n}$ $\in H^{2}(M;\mathbb{R})$ generate the cohomology ring of $M$ and $\omega$ is written as $\omega=\sum\limits_{i=1}^{n}t_{i}w_{i}$, where $t_{i} \in \mathbb{R}$ for $1\leq i \leq n$. Then the cohomology ring $H^{*}(M;\mathbb{R})$ is isomorphic to the Poincar$\acute{{\text e}}$ duality algebra $\mathbb{R}[\partial_{1},\ldots, \partial_{n}]/{\rm Ann}({\rm vol}(M))$. 
\end{lemma}

Considering the special case of Lemma 2.5, we have the following result.
\begin{lemma}
Under the assumptions of Lemma 2.5, assume that the symplectic volume of $M$ is a polynomial of $t_{i_{1}}+\cdots+t_{i_{k}}\ (1< k < n)$ and $t_{j_{a}}\  (j_{a}\neq i_{1}, \ldots ,i_{k}, 1 \leq a \leq n-k)$ for some $k$.
Then we have the relation $w_{i_{l}}=w_{i_{m}}$ in $H^{*}(M;\mathbb{R})$ for all $1\leq l,m \leq k$.
\end{lemma}

\begin{proof}
From the assumption, we have
\[
\int_{M}e^{\omega}=f(t_{i_{1}}+\cdots+t_{i_{k}},t_{j_{1}},\ldots, t_{j_{n-k}}).
\]
Since  $\omega=\sum\limits_{i=1}^{n}t_{i}w_{i}$, we have, for $1 \leq l$, $m \leq k$, 
\begin{eqnarray*}
\int_{M}(w_{i_{l}}-w_{i_{m}})e^{\omega}&=&\left(\frac{\partial}{\partial t_{i_{l}}}-\frac{\partial}{\partial t_{i_{m}}}\right)\int_{M}e^{\omega}\\
&=&\left(\frac{\partial}{\partial t_{i_{l}}}-\frac{\partial}{\partial t_{i_{m}}}\right)f(t_{i_{1}}+\cdots+t_{i_{k}},t_{j_{1}},\ldots, t_{j_{n-k}})\\
&=&0.
\end{eqnarray*}
\end{proof}

\section{Applications to the topology of weight varieties}

\subsection{Multiple weight varieties}
In this subsection, we give the definition of multiple weight varieties. 
Let $G$ be a compact, connected Lie group with Lie algebra $\mathfrak{g}$, and let $T$ be a maximal torus of $G$ with Lie algebra $\mathfrak{t}$. Let $l={\rm dim}\ T$ and let $\mathfrak{g}^{*}$ and $\mathfrak{t}^{*}$ be the dual vector spaces of $\mathfrak{g}$ and $\mathfrak{t}$ respectively. We identify $\mathfrak{g}$ and $\mathfrak{t}$ with $\mathfrak{g}^{*}$ and $\mathfrak{t}^{*}$ by an invariant inner product, respectively.

Let $\Delta$ be the root system of $\mathfrak{g}$. Let $\Delta_{+}$ be a set of positive roots, let $\{\alpha_{1},  \ldots, \alpha_{l} \}\subset \Delta_{+}$ be the set of simple roots, let $\{ \Lambda_{1},\ldots, \Lambda_{l}\}$ be the set of fundamental weights.

Let us set
\[
\mathfrak{t}^{*}_{+}:= \mathbb{R}_{\geq 0}\Lambda_{1}+\mathbb{R}_{\geq 0}\Lambda_{2}+\cdots +\mathbb{R}_{\geq 0}\Lambda_{l}, \hspace{5pt}
\mathfrak{t}^{*}_{++}:= \mathbb{R}_{> 0}\Lambda_{1}+\mathbb{R}_{> 0}\Lambda_{2}+\cdots +\mathbb{R}_{> 0}\Lambda_{l},
\]
\[
P_{+}:= \mathbb{Z}_{\geq 0}\Lambda_{1}+\mathbb{Z}_{\geq 0}\Lambda_{2}+\cdots +\mathbb{Z}_{\geq 0}\Lambda_{l}, \hspace{5pt}
P_{++}:= \mathbb{Z}_{> 0}\Lambda_{1}+\mathbb{Z}_{> 0}\Lambda_{2}+\cdots +\mathbb{Z}_{> 0}\Lambda_{l}.\\
\]

The set $\mathfrak{t}^{*}_{+}$ is called the Weyl chamber and it is a fundamental domain of action of the Weyl group $W$ on $\mathfrak{t}^{*}$.  Elements in $P_{+}$ are called dominant weights.
The left coadjoint action of $G$ on $\mathfrak{g}^{*}$ is defined by $g \cdot f={\rm Ad}^{*}(g)f$  for all $g \in G$ and for all $f \in \mathfrak{g}^{*}$, where 
\[
\langle {\rm Ad}^{*}(g)f, X \rangle := \langle f, {\rm Ad}(g^{-1})X \rangle \hspace{5pt}  (X \in \mathfrak{g}).
\]

 Let ${\mathcal O}_{\lambda}$ be the coadjoint orbit through $\lambda \in \mathfrak{t}^{*}_{+}$. Then the intersection ${\mathcal O}_{\lambda} \cap \mathfrak{t}^{*}$ is the  $W$-orbit of $\lambda$ and the set ${\mathcal O}_{\lambda}\cap \mathfrak{t}^{*}_{+}$ is equal to $\{\lambda \}$. In other words, coadjoint orbits are parametrized by elements in $\mathfrak{t}^{*}_{+}$. We denote by $G_{\lambda}$ the isotropy subgroup of $\lambda \in \mathfrak{t}_{+}^{*}$ for the coadjoint action of $G$ on $\mathfrak{g}^{*}$.
 
Each coadjoint orbit ${\mathcal O}_{\lambda}$ has a $G$-invariant symplectic structure $\omega_{\lambda}$ called the Kostant-Kirillov-Souriau symplectic form defined by
\[
(\omega_{\lambda})_{x}(\xi_{X},\xi_{Y})=\langle x, [X,Y] \rangle \hspace{2em} {\rm for} \hspace{1em} x \in {\mathcal O}_{\lambda} \hspace{1em} {\rm and} \hspace{1em} X,Y \in \mathfrak{g}.
\]
The action of $G$ on ${\mathcal O}_{\lambda}$ is Hamiltonian, and the moment map $\iota_{\lambda}:{\mathcal O}_{\lambda} \to \mathfrak{g}^{*}$ is given by the inclusion.
The action of the maximal torus $T$ of $G$ on the coadjoint orbit ${\mathcal O}_{\lambda}$ is also Hamiltonian, and the moment map $\Phi_{\lambda}: {\mathcal O}_{\lambda} \to \mathfrak{t}^{*}$ is given by the composition of the inclusion $\iota_{\lambda}: {\mathcal O}_{\lambda} \to \mathfrak{g}^{*}$ and the projection $\pi:\mathfrak{g}^{*} \to \mathfrak{t}^{*}$.

More generally, for $\lambda_{1}, \lambda_{2}, \ldots,\lambda_{n} \in \mathfrak{t}_{+}^{*}$, the diagonal action of $T$ on the product of the coadjoint orbits ${\mathcal O}_{\lambda_{1}} \times \cdots \times {\mathcal O}_{\lambda_{n}}$ is also Hamiltonian and the moment map $\Phi:{\mathcal O}_{\lambda_{1}} \times \cdots \times {\mathcal O}_{\lambda_{n}} \to \mathfrak{t}^{*}$ is given by  $\Phi(x_{1},\ldots ,x_{n}) ={\displaystyle \sum_{i=1}^{n}} \Phi_{\lambda_{i}}(x_{i})$.

For $\mu \in \mathfrak{t}^{*}$, a multiple weight variety is defined by
\[
{\mathcal M}_{T}:=\Phi^{-1}(\mu)/T=\left\{(x_{1},\ldots ,x_{n}) \in {\mathcal O}_{\lambda_{1}}\times \cdots \times {\mathcal O}_{\lambda_{n}}\ \middle| \  \sum_{i=1}^{n}\Phi_{\lambda_{i}}(x_{i})=\mu \right\}/T.
\]

We consider the case $G=SU(l+1)$ and $\lambda_{1},\ldots, \lambda_{n} \in \mathfrak{t}^{*}_{++}$. The action of the center $Z(G)$ on $\mathcal{O}_{\lambda_{1}}\times \cdots \times \mathcal{O}_{\lambda_{n}}$ is trivial. Considering $T/Z(G)$ instead of $T$, for a regular value $\mu \in \mathfrak{t}^{*}$ of the moment map, there exists an open neighborhood $U$ of $\mu$ such that $T/Z(G)$ acts freely on $\Phi^{-1}(U)$. Thus, $\mathcal{M}_{T}$ is a compact smooth symplectic manifold in this case.

As a special case of Theorem 2.1 and Corollary 2.2, we have the following.

\begin{theorem}
Suppose that $G=SU(l+1)$, $\lambda_{i}=\sum\limits_{j=1}^{l}p_{i,j}\alpha_{j} \in \mathfrak{t}^{*}_{++}$ $(1\leq i \leq n,\ p_{i,j} \in \mathbb{R})$ and $\mu_{0} \in \mathfrak{t}^{*}$ is a regular value of the moment map $\Phi:\mathcal{O}_{\lambda_{1}}\times \cdots \times \mathcal{O}_{\lambda_{n}}\to \mathfrak{t}^{*}$. Let $\omega_{\mu_{0}}$ be the symplectic form on the smooth symplectic manifold $\mathcal{M}_{T}=\Phi^{-1}(\mu_{0})/T$. Then there exist an open neighborhood $U$ of $\mu_{0}$ and a diffeomorphism $\phi_{\mu}:\Phi^{-1}(\mu_{0})/T \to \Phi^{-1}(\mu)/T$ for all $\mu \in U$. Moreover there exist cohomology classes $z_{i,j}$ and $v_{j} \in H^{2}(\mathcal{M}_{T};\mathbb{R})$ that generate the cohomology ring $H^{*}(\mathcal{M}_{T}; \mathbb{R})$ multiplicatively such that 
\[
\phi^{*}_{\mu}([\omega_{\mu}])=\sum_{i=1}^{n}\sum_{j=1}^{l}p_{i,j}z_{i,j}+\sum_{j=1}^{l}x_{j}v_{j}
\]
for all $\mu=\sum\limits_{j=1}^{l}x_{j}\alpha_{j} \in U \subset \mathfrak{t}^{*}$, where $x_{j}\in \mathbb{R}$.
\end{theorem}

\begin{proof}
The symplectic form on ${\mathcal O}_{\lambda_{1}}\times \cdots \times {\mathcal O}_{\lambda_{n}}$ is given by $\omega=\sum\limits_{i=1}^{n}{\rm pr^{*}_{i} }\omega_{\lambda_{i}}$ where $\omega_{\lambda_{i}}$ is the symplectic form on ${\mathcal O}_{\lambda_{i}}$ and ${\rm pr}_{i}$ $(1\leq i \leq n)$ is the $i$-th projection ${\rm pr}_{i}:{\mathcal O}_{\lambda_{1}}\times \cdots \times {\mathcal O}_{\lambda_{n}} \to {\mathcal O}_{\lambda_{i}}$.

When $\lambda \in \mathfrak{t}^{*}_{++}$,  then $G_{\lambda}=T$. Since $G_{\lambda_{i}} \subset G_{\alpha_{j}}$, the $G$-invariant closed 2-form $\omega_{i,j}$ defined by
\[
(\omega_{i,j})_{g \cdot \lambda_{i}}(\xi_{X},\xi_{Y})=\langle g\cdot \alpha_{j}, [X,Y]\rangle \hspace{15pt}X,Y \in \mathfrak{g},\  g\in G, \ g\cdot \lambda_{i} \in {\mathcal O}_{\lambda_{i}}
\]
is well-defined for all $1\leq i \leq n$ and $1\leq j\leq l$.
Similarly, the map $\Phi_{i,j}:\mathcal{O}_{\lambda_{i}}\to \mathfrak{t}^{*}$ defined by
\[
\Phi_{i,j}(g\cdot \lambda_{i})=\pi(g\cdot \alpha_{j})
\]
is also well-defined for all $1\leq i \leq n$ and $1\leq j\leq l$.
Since $\lambda_{i}=\sum\limits_{j=1}^{l}p_{i,j}\alpha_{j} \in \mathfrak{t}^{*}_{++}$ $(1\leq i \leq n)$, the symplectic form $\omega=\sum\limits_{i=1}^{n}{\rm pr^{*}_{i} }\omega_{\lambda_{i}}$ and the moment map $\Phi=\sum\limits_{i=i}^{n}\Phi_{\lambda_{i}}$ are written as
\[
\omega=\sum\limits_{j=1}^{l}{\rm pr^{*}_{i}}(p_{1,j}\omega_{1,j})+\cdots+\sum\limits_{j=1}^{l}{\rm pr^{*}_{n}}(p_{n,j}\omega_{n,j}),
\]
\[
\Phi=\sum\limits_{j=1}^{l}p_{1,j}\Phi_{1,j}+\cdots+\sum\limits_{j=1}^{l}p_{n,j}\Phi_{n,j}.
\]
Then, $\omega_{i,j}$ and $\Phi_{i,j}$ satisfy the condition
\[
\iota_{\xi_{X}}\omega_{i,j}=d\langle \Phi_{i,j}, X \rangle \hspace{15pt} X \in \mathfrak{t}
\]
for all $1\leq i \leq n$ and $1\leq j \leq l$.
Since $\lambda_{i}\in \mathfrak{t}^{*}_{++}$ $(1\leq i\leq n)$, the coadjoint orbit $\mathcal{O}_{\lambda_{i}}$ is diffeomorphic to the flag variety $G/T$. The equivariant cohomology $H^{*}_{T}(\mathcal{O}_{\lambda_{i}};\mathbb{R})$ is generated by the equivariant cohomology classes $[\omega_{i,j}-\Phi_{i,j}]$ $(1\leq j \leq l)$ as a  $H^{*}(BT;\mathbb{R})$-module ([18]). Using the K$\ddot{{\rm u}}$nneth formula ([13], [8]), we see that $ H^{*}_{T}(\mathcal{O}_{\lambda_{1}}\times \cdots \times \mathcal{O}_{\lambda_{n}};\mathbb{R})$ is also generated by $[\omega_{i,j}-\Phi_{i,j}]\ (1\leq i\leq n,\ 1\leq l \leq l)$ as a $H^{*}(BT;\mathbb{R})$-module. Thus, $\mathcal{O}_{\lambda_{1}}\times \cdots \times \mathcal{O}_{\lambda_{n}}$ satisfies the assumptions in Theorem 2.1 and Corollary 2.2.

\end{proof}
By Theorem 3.1, we see ${\rm vol}(\mathcal{M}_{T})$ is a polynomial of $p_{i,j}$ and $x_{j}$ $(1\leq i \leq n,\ 1\leq j\leq l).$
\subsection{Flow polytopes and volumes}
In this subsection we review the definition of the flow polytope associated to the root system of type $A$ based on [3] and [14].

Let $e_{1},\ldots, e_{l+1}$ be the standard basis of $\mathbb{R}^{l+1}$,
\[
A_{l}^{+}=\{e_{i}-e_{j}\ |\  1\leq i < j \leq l+1\}
\]
be the positive root system of type $A$ with rank $l$ and
\[
\alpha_{1}=e_{1}-e_{2},\  \alpha_{2}=e_{2}-e_{3}, \ldots,\  \alpha_{l}=e_{l}-e_{l+1}
\]
be the simple roots. We assign a positive integer $m_{i,j}$ to each pair $(i,j)$ and set $m=(m_{i,j})$ and $M=\sum\limits_{1\leq i<j\leq l+1} m_{i,j}$.

\begin{definition}
Let $C(A_{l}^{+})$ be the convex cone generated by $A_{l}^{+}$:
\[
C(A_{l}^{+})=\{q_{1}\alpha_{1}+\cdots +q_{l}\alpha_{l}\ |\ q_{1},\ldots, q_{l}\in \mathbb{R}_{\geq0}\}.
\]
For $h \in C(A_{l}^{+})$, the polytope
\[
P_{l,m}(h)=
\left\{ 
(x_{i,j}^{k})\in \mathbb{R}^{M}
\ \middle| \ 
\begin{aligned}
1\leq i < j \leq l+1,\ \  1 \leq k \leq m_{i,j}, \\
 x_{i,j}^{k} \geq 0,\ \   \sum_{i<j}\sum_{k}x_{i,j}^{k}(e_{i}-e_{j})=h
\end{aligned}
\right\}
\]
is called the flow polytope.
\end{definition}

For $h \in C(A_{l}^{+})\cap \sum\limits_{i=1}^{l}\mathbb{Z}\alpha_{i}$, the Kostant partition function $p_{l,m}$ is defined by
\[
p_{l,m}(h)=|P_{l,m}(h) \cap \mathbb{Z}^{M}|
\]
and the function $h \mapsto p_{l,m}(h)$ becomes a polynomial ([3]).

The volume $v_{l,m}(h)$ of a flow polytope $P_{l,m}(h)$ is a piecewise polynomial on $C(A_{l}^{+})$. In particular, $v_{l,m}(h)$ is given by
\[
v_{l,m}(h)=\lim_{k\to \infty}\frac{1}{k^{d}}\cdot p_{l,m}(kh)
\]
for $h \in C(A_{l}^{+})\cap \sum\limits_{i=1}^{l}\mathbb{Z}\alpha_{i}$, where $d$ is the degree of $p_{l,m}$.

Let $n$ be a positive integer. When $m_{i,j}=n$ for all $i, j$, we denote the flow polytope and its volume by $P_{l,n}(h)$ and $v_{l,n}(h)$ instead of $P_{l,m}(h)$ and $v_{l,m}(h)$ respectively. 
  
\subsection{Symplectic volume of special weight varieties}
Let us set $\rho=\frac{1}{2}\sum\limits_{i=1}^{l}\alpha_{i}$.
For $\lambda_{1},\ldots, \lambda_{n}\in \mathfrak{t}^{*}_{+}-\{0\}$ and $\mu \in \mathfrak{t}^{*}$, we say that $\mu$ is sufficiently close to $\lambda=\lambda_{1}+\cdots+\lambda_{n}$ when the following condition holds:

If $\sigma_{1}(\lambda_{1}+\rho)+\cdots+\sigma_{n}(\lambda_{n}+\rho)-(\mu+n\rho) \in C(A_{l}^{+})$ for $\sigma_{1},\ldots, \sigma_{n}\in W$, then $\sigma_{1}=\cdots=\sigma_{n}=e$.

\begin{theorem}
Suppose that $G=SU(l+1)$, $\lambda_{1}, \ldots, \lambda_{n} \in \mathfrak{t}_{+}^{*}-\{0\}$ and let $\mu \in \mathfrak{t}^{*}$ be a a regular value of the moment map sufficiently close to $\lambda=\lambda_{1}+\cdots +\lambda_{n}$. Then
\[
{\rm vol}(\mathcal{M}_{T})=v_{l,n}(\lambda-\mu).
\]
\end{theorem}

\begin{proof}
First, we consider the case $\lambda_{1}, \ldots, \lambda_{n}\in P_{+}-\{0\}$ and $\mu \in P$. Recall that a multiple weight variety is defined by
\[
\mathcal{M}_{T}(\lambda_{1},\ldots, \lambda_{n},\mu):=\Phi^{-1}(\mu)/T
\]
As proved in [15], the symplectic volume of a multiple weight variety is given by the asymptotic behavior of the weight multiplicity, that is, 
\begin{align}
{\rm vol}(\mathcal{M}_{T}(\lambda_{1},\ldots, \lambda_{n},\mu))=\lim_{k \to \infty}\frac{1}{k^{d}}\cdot[V_{k\lambda_{1}}\otimes \cdots \otimes V_{k\lambda_{n}};W_{k\mu}]
\end{align}
where $V_{\lambda}$ is the irreducible representation of $G$ with highest weight $\lambda$, $W_{\mu}$ is the weight space associated with $\mu$ and $d$ is the dimension of $\mathcal{M}_{T}$.

On the other hand, a generalization of the Kostant multiplicity formula ([7]) gives us
\begin{align*}
&[V_{\lambda_{1}}\otimes \cdots \otimes V_{\lambda_{n}};W_{\mu}]\\
&=\sum\limits_{\sigma_{1},\ldots, \sigma_{n}\in W}\left(\prod\limits_{i=1}^{n}\varepsilon(\sigma_{i})\right) p_{l,n}(\sigma_{1}(\lambda_{1}+\rho)+\cdots+\sigma_{n}(\lambda_{n}+\rho)-(\mu+n\rho))
\end{align*}
where $\varepsilon(\sigma_{i})$ is the signature of $\sigma_{i} \in$ W.

When $\mu$ is sufficiently close to $\lambda$, we have 
\[
\sum\limits_{\sigma_{1},\ldots, \sigma_{n}\in W}\left(\prod\limits_{i=1}^{n}\varepsilon(\sigma_{i})\right) p_{l,n}(\sigma_{1}(\lambda_{1}+\rho)+\cdots+\sigma_{n}(\lambda_{n}+\rho)-(\mu+n\rho))=p_{l,n}(\lambda-\mu).
\]
Therefore, we have
\[
{\rm vol}(\mathcal{M}_{T}(\lambda_{1},\ldots, \lambda_{n},\mu))
=
\lim_{k\to \infty}\frac{1}{k^{d}}\cdot p_{l,n}(k(\lambda-\mu))
=v_{l,n}(\lambda-\mu).
\]

Next, we consider the case $\lambda_{1},\ldots, \lambda_{n}\in (P\otimes \mathbb{Q})\cap \mathfrak{t}^{*}_{+}-\{0\}$ and $\mu \in P\otimes \mathbb{Q}$. Then, there exists $m\in \mathbb{Z}_{>0}$ such that $m\lambda_{i}\in P_{+}-\{0\}$ and $m\mu \in P$. Moreover, $m\mu$ is sufficiently close to $m\lambda$. Thus we have
\[
{\rm vol}(\mathcal{M}_{T}(m\lambda_{1},\ldots,m\lambda_{n}, m\mu))=\lim_{k\to \infty}\frac{1}{k^{d}}\cdot p_{l,n}(km(\lambda-\mu)).
\]
Since 
\[
{\rm vol}(\mathcal{M}_{T}(m\lambda_{1},\ldots,m\lambda_{n}, m\mu))=m^{d}{\rm vol}(\mathcal{M}_{T}(\lambda_{1},\ldots,\lambda_{n}, \mu))
\] and 
\[
\lim_{k\to \infty}\frac{1}{k^{d}}\cdot p_{l,n}(km(\lambda-\mu))=m^{d}v_{l,n}(\lambda-\mu),
\]
we have
\[
{\rm vol}(\mathcal{M}_{T}(\lambda_{1},\ldots, \lambda_{n},\mu))=v_{l,n}(\lambda-\mu)
\]
for $\lambda_{1},\ldots, \lambda_{n}\in (P\otimes \mathbb{Q})\cap \mathfrak{t}^{*}_{+}-\{0\}$ and $\mu \in P\otimes \mathbb{Q}$.

By Theorem 3.1, ${\rm vol}(\mathcal{M}_{T}(\lambda_{1},\ldots, \lambda_{n},\mu))$ and $v_{l,n}(\lambda-\mu)$ are polynomial. By continuity of volume functions, above equation holds for $\lambda_{1}, \ldots, \lambda_{n} \in \mathfrak{t}_{+}^{*}-\{0\}$ and $\mu \in \mathfrak{t}^{*}$.
\end{proof}

\begin{example}
We consider the case $G=SU(3)$ and  $n=2$. In [15], a formula for ${\rm vol}(\mathcal{M}_{T})$ is obtained by calculating the right hand side of (3.1) directly.
For example, $\lambda_{1}=p\alpha_{1}+q\alpha_{2}$ ,$\lambda_{2}=r\alpha_{1}+s\alpha_{2} \in \mathfrak{t}^{*}_{++}$ and $\mu=x\alpha_{1}+y\alpha_{2} \in \mathfrak{t}^{*}$ is in the alcove which is adjacent to the $\lambda_{1}+\lambda_{2}$, where $p,q,r,s,x,y \in \mathbb{R}$. Then, ${\rm vol}(\mathcal{M}_{T})$ is given by
\[
{\rm vol}(\mathcal{M}_{T})=
\begin{cases}
\frac{1}{12}(q+s-y)^{3}(2(p+r-x)-(q+s-y))\\ \hspace{80pt} ({\rm if}\  \lambda_{1}+\lambda_{2}-\mu \in \mathbb{R}_{>0}\alpha_{1}+\mathbb{R}_{>0}(\alpha_{1}+\alpha_{2}))\\
\frac{1}{12}(p+r-x)^{3}(-(p+r-x)+2(q+s-y))\\ \hspace{80pt} ({\rm if}\  \lambda_{1}+\lambda_{2}-\mu \in \mathbb{R}_{>0}(\alpha_{1}+\alpha_{2})+\mathbb{R}_{>0}\alpha_{2}).\\
\end{cases}
\]

We can calculate intersection pairings  using Corollary 2.3. By Theorem 3.1, we have the cohomology classes $z_{1,1}, z_{1,2}, z_{2,1}, z_{2,2}, v_{1}, v_{2} \in H^{2}(\mathcal{M}_{T}; \mathbb{R})$. Using the fomula (2.2), we have
\begin{eqnarray*}
{\rm vol}(\mathcal{M}_{T})&=&\frac{1}{4!}\int_{\mathcal{M}_{T}}(pz_{1,1}+qz_{1,2}+rz_{2,1}+sz_{2,2}+xv_{1}+yv_{2})^{4}\\
&=&\sum\frac{p^{d_{1}}\cdots y^{d_{6}}}{d_{1}!\cdots d_{6}!}\int_{\mathcal{M}_{T}}z_{1,1}^{d_{1}}\cdots v_{2}^{d_{6}}
\end{eqnarray*}
where $d_{1},\ldots, d_{6}\in \mathbb{Z}_{\geq0}$ and $\sum\limits_{i=1}^{6}d_{i}=\frac{1}{2}{\rm dim}_{\mathbb{R}}\mathcal{M}_{T}=4$. Note also that
\[
\int_{\mathcal{M}_{T}}z_{1,1}^{d_{1}}\cdots v_{2}^{d_{6}}=\frac{\partial^{4}}{\partial p^{d_{1}}\cdots \partial y^{d_{6}}}{\rm vol}(\mathcal{M}_{T}).
\]
We consider the case 
$ \lambda_{1}+\lambda_{2}-\mu \in \mathbb{R}_{>0}(\alpha_{1}+\alpha_{2})+\mathbb{R}_{>0}\alpha_{2}$.
Then, we have 
\begin{align*}
&\int_{\mathcal{M}_{T}}z_{1,1}^{4}=-2, \hspace{2pt}\int_{\mathcal{M}_{T}}z_{1,1}^{3}z_{1,2}=1,
\hspace{2pt} \int_{\mathcal{M}_{T}}z_{1,1}^{3}z_{2,1}=-2,\\
&\int_{\mathcal{M}_{T}}z_{1,1}^{3}z_{2,2}=1, \hspace{2pt} \int_{\mathcal{M}_{T}}z_{1,1}^{3}v_{1}=2,  
\hspace{2pt}\int_{\mathcal{M}_{T}}z_{1,1}^{3}v_{2}=-1
\end{align*}
and so on.

In particular, we can obtain the Betti numbers using the above calculations.
Since $\mathcal{M}_{T}$ is a compact connected symplectic manifold of dimension eight and $z_{1,1}, z_{1,2},$ $z_{2,1}, z_{2,2}, v_{1}, v_{2} \in H^{2}(\mathcal{M}_{T}; \mathbb{R})$ generate $H^{*}(\mathcal{M}_{T}; \mathbb{R})$ multiplicatively, we have dim $H^{0}(\mathcal{M}_{T};\mathbb{R})={\rm dim}\ H^{8}(\mathcal{M}_{T};\mathbb{R})=1$ and dim $H^{odd}(\mathcal{M}_{T};\mathbb{R})=0$.

To calculate ${\rm dim}\ H^{2}(\mathcal{M}_{T},\mathbb{R})$, we consider the equation
\[
a_{1}z_{1,1}+a_{2}z_{1,2}+a_{3}z_{2,1}+a_{4}z_{2,2}+a_{5}v_{1}+a_{6}v_{2}=0,
\]
where $a_{i}\in \mathbb{R} \ \  (0\leq i \leq 6)$.
Multiplying by $z_{1,1}^{3}$ and integrating both sides, we have 
\[
\int_{\mathcal{M}_{T}}(a_{1}z_{1,1}^{4}+a_{2}z_{1,1}^{3}z_{1,2}+a_{3}z_{1,1}^{3}z_{2,1}+a_{4}z_{1,1}^{3}z_{2,2}+a_{5}z_{1,1}^{3}v_{1}+a_{6}z_{1,1}^{3}v_{2})=0.
\]
Using the above computations, we have 
\[
-2a_{1}+a_{2}-2a_{3}+a_{4}+2a_{5}-a_{6}=0.
\] 
Similarly, multiplying by $z_{1,1}z_{1,2}z_{2,1}$ and integrating both sides, we have
\[
a_{1}+a_{3}-a_{5}=0.
\]
Multiplying by any monomial of $z_{1,1}, z_{1,2}, z_{2,1}, z_{2,2}, v_{1}, v_{2}$ of degree three  and integrating both sides, we always obtain constant multiples of the above two equations.
Therefore, we have
\[
{\rm dim}\ H^{2}(\mathcal{M}_{T},\mathbb{R})={\rm rank}
\left(
\begin{array}{cccccc}
-2&1&-2&1&2&-1\\
1&0&1&0&-1&0\\
\end{array}
\right)=2.
\]
Analogously, we have ${\rm dim}\ H^{4}(\mathcal{M}_{T},\mathbb{R})=2$ by simular argument as above.

Therefore the Poincar$\acute{{\text e}}$ polynomial $P_{t}({\mathcal M}_{T})$ is equal to $1+2t^{2}+2t^{4}+2t^{6}+t^{8}$. 
\end{example}

\subsection{Cohomology rings of special weight varieties}
In this subsection, we consider special flow polytopes $P_{l,m}(h)$ based on [14], where $h\in C(A_{l}^{+})$ is in the nice chamber (see Definition 3.5 below). In this case, the volume function $v_{l,m}(h)$ is characterized by a system of differential equations. Applying this characterization of $v_{l,m}(h)$, we have an explicit presentation of the cohomology ring of the multiple weight variety $\mathcal{M}_{T}$ of special type.

\begin{definition}[{[3]}]
The open subset $\mathfrak{c}_{{\rm nice}}$ of $C(A_{l}^{+})$ is defined by
\[
\mathfrak{c}_{{\rm nice}}=\left\{h:=q_{1}\alpha_{1}+\cdots+q_{l}\alpha_{l} \middle| \  q_{1},\ldots , q_{l} \in \mathbb{R}_{\geq 0},\  0<q_{1}<q_{2}<\cdots <q_{l} \right\},
\]
which is called the nice chamber.
\end{definition}

In [14], a characterization of the volume function $v_{l,m}(h)$ is given in the special case that $h \in C(A_{l}^{+})$ is in the nice chamber. 

\begin{theorem}[{[14]}]
If $h=q_{1}\alpha_{1}+\cdots +q_{l}\alpha_{l} \in C(A_{l}^{+})$ is in the nice chamber, then $v=v_{l,m}(h)$ satisfies the following system of differential equations
\[
\left\{
\begin{array}{l}
\partial^{m_{l,l+1}}_{l}v=0\\
\partial^{m_{l-1,l}}_{l-1}(\partial_{l-1}+\partial_{l})^{m_{l-1,l+1}}v=0\\
\vdots\\
\partial_{1}^{m_{1,2}}(\partial_{1}+\partial_{2})^{m_{1,3}}\cdots (\partial_{1}+\cdots +\partial_{l})^{m_{1,l+1}}v=0,
\end{array}
\right.
\]
where $\partial_{i}=\frac{\partial}{\partial q_{i}}$ for $1\leq i \leq l$. Conversely, the polynomial of degree $M-l$ satisfying the above equations is equal to a constant multiple of $v_{l,m}(h)$.
\end{theorem}

\begin{remark}
The coordinates in [14] and [3] are $(a_{1},\ldots, a_{l})$ but we consider  the coordinates $(q_{1},\ldots, q_{l})$ where $q_{1}=a_{1},\ q_{2}=a_{1}+a_{2},\ldots,\ q_{l}=a_{1}+a_{2}+\cdots+a_{l}$.
\end{remark}

\begin{theorem}
If $h=q_{1}\alpha_{1}+\cdots+q_{l}\alpha_{l} \in C(A_{l}^{+})$ is in the nice chamber, then ${\rm Ann}(v)$ is generated by the differential operators in Theorem 3.6.
\end{theorem}

\begin{proof}
Let $M_{i}={\displaystyle \sum_{j=i+1}^{l+1}m_{i,j}}$ and let $J$ be the ideal in $\mathbb{R}[\partial_{1},\ldots,\partial_{l}]$ generated by differential operators in Theorem 3.6. By Theorem 3.6, for each $i\in \{1,\ldots, l\}$ there exists a homogeneous polynomial $P_{i}(\partial_{i},\ldots, \partial_{l})$ of degree $M_{i}$ such that 
\begin{align}
\partial_{i}^{M_{i}}v=P_{i}(\partial_{i},\ldots, \partial_{l})v,
\end{align} 
and such that the degree of $\partial_{i}$ in $P_{i}(\partial_{i},\ldots, \partial_{l})$ is less than $M_{i}$.
Note that $\partial_{i}^{M_{i}}-P_{i}(\partial_{i},\ldots , \partial{_{l}})$ belongs to $J$.

Let $D$ be an operator in ${\rm Ann}(v)$ and let us write 
\[
D=\sum_{c\geq 0}D_c
\]
where $D_c$ is the homogeneous part of degree $c$ in $D$. 
Since $v$ is homogeneous of degree $M-l$, we see that 
$D_cv$ is homogeneous of degree $M-l-c$ for each $c$ with $0\leq c\leq M-l$, 
while $D_cv=0$ for $c>M-l$. 
Hence $Dv=0$ implies $D_cv=0$ for all $c$. 
Therefore, we can assume, without loss of generality, that $D=D_{c}$ for some $c$.
Moreover, using (3.2), we have
\[
D\equiv D'\ \ \     {\rm mod}\ J
\]
and
\begin{align}
D':=\sum_{ |\alpha|=c}\beta_{\alpha}\partial^{\alpha}
\end{align}
where $\partial=(\partial_{1},\ldots,\partial_{l})$, $\beta_{\alpha} \in \mathbb{R}$,  $\alpha=(\alpha_{1},\ldots, \alpha_{l}) \in \mathbb{Z}^{l}_{\geq 0}$ is a multi-index satisfying $\alpha_{i}<M_{i}$ for all $i\in\{1,\ldots, l\}$ and $|\alpha|:=\alpha_{1}+\cdots+\alpha_{l}$. 

We want to show that $D'=0$. Suppose, on the contrary, that $D' \neq 0$, so that
\[
A:=\left\{\alpha \  \middle| \ \beta_{\alpha}\neq0 \right\}
\]
in (3.3) is not empty. 
Taking the lexicographic maximum $I=(i_{1},\ldots,i_{l})$ of the multi-indices $\alpha$ in A, we can write
\[
D'v=\left(\partial^{I}+\sum_{\alpha \prec I }\beta_{\alpha} \partial^{\alpha}\right)v=0,
\]
(we can assume, without loss of generality that the coefficient $\beta_{I}$ in $D'$ is 1).
Then we have

\begin{align}
\partial_{1}^{M_{1}-1}\cdots \partial_{l}^{M_{l}-1}v&=\partial_{1}^{M_{1}-1-i_{1}}\cdots \partial_{l}^{M_{l}-1-i_{l}}\partial^
{I}v\\
&=-\left(\partial_{1}^{M_{1}-1-i_{1}}\cdots \partial_{l}^{M_{l}-1-i_{l}}\sum_{\alpha \prec I}\beta_{\alpha}\partial^{\alpha} \right)v \notag \\
&=-\left(\sum_{\alpha \prec I}\beta_{\alpha}\partial_{1}^{M_{1}-1-i_{1}+\alpha_{1}}\cdots \partial_{l}^{M_{l}-1-i_{l}+\alpha_{l}}\right)v \notag.
\end{align}

From  [14, Proposition 17], we know that the coefficient of $q_{1}^{M_{1}-1}\cdots q_{l}^{M_{l}-1}$ in the volume function $v$ is nonzero and so the left hand side of (3.4) is nonzero.

On the other hand, the right hand side is zero by the claim below. Therefore, we obtain a contradiction.
\end{proof}

\vspace{10pt}
\hspace{-18.5pt}
{\bf Claim.} {\it If $\alpha=(\alpha_{1},\ldots, \alpha_{l}) \in \mathbb{Z}_{\geq 0}^{l}$ satisfies $|\alpha|=M_{1}+\cdots+M_{l}-l$ and $\alpha \prec (M_{1}-1,\ldots ,M_{l}-1)$, then $\partial^{\alpha}v=0$. }
\vspace{10pt}

\begin{proof}
Since $\alpha \prec (M_{1}-1,\ldots ,M_{l}-1)$, there exists $s_{\alpha} \in \{1,\ldots, l\}$ such that $\alpha_{s_{\alpha}}<M_{s_{\alpha}}-1$ and $\alpha_{i}=M_{i}-1$ for all $i<s_{\alpha}$. Since $|\alpha|=M_{1}+\cdots+M_{l}-l$, there exists $j \in \{s_{\alpha}+1,\ldots, l\}$ such that $\alpha_{j}>M_{j}-1$ and $\alpha_{j'}\leq M_{j'}-1$ for all $j'<j$.

Using (3.2), $\partial_{j}^{\alpha_{j}}$ (mod $J$) can be written as a homogeneous polynomial of $\partial_{j},\ldots, \partial_{l}$ such that the degree of $\partial_{j}$ is less than $M_{j}$.
Therefore, we can write
\[
\partial^{\alpha}\equiv \sum_{\beta}\gamma_{\beta}\partial_{1}^{\alpha_{1}}\cdots \partial_{j-1}^{\alpha_{j-1}}\partial_{j}^{\beta_{j}}\cdots \partial_{l}^{\beta_{l}}\ \ \ \ {\rm mod}\ J
\]
where $\gamma_{\beta} \in \mathbb{R}$ and $\beta=(\beta_{j},\ldots,\beta_{l})\in \mathbb{Z}_{\geq 0}^{l-j+1}$is a multi-index satisfying $\beta_{j}\leq M_{j}-1$.

Since $\alpha_{1}+\cdots+\alpha_{j-1}+\beta_{j}+\cdots+\beta_{l}=M_{1}+\cdots+M_{l}-l$, $\alpha_{s_{\alpha}}<M_{s_{\alpha}}-1$,  $\alpha_{i'}\leq M_{i'}-1$ for all $i' \in \{1,\ldots,j-1\}$ and $\beta_{j}\leq M_{j}-1$, there exists  $k \in \{j+1,\ldots, l\}$ such that $\beta_{k}>M_{k}-1$ and $\beta_{k'}\leq M_{k'}-1$ for all $k'<k$. 

Repeating this argument, $\partial^{\alpha}$ (mod $J$)  can be written as a homogeneous polynomial of $\partial_{1},\ldots,\partial_{l}$ such that the degree of $\partial_{l}$ is greater than $M_{l}-1$. Since $\partial_{l}^{M_{l}}v=0$, we have $\partial^{\alpha}v=0$.
\end{proof}

The case all $m_{i,j}=n$ $(1\leq i<j\leq l+1)$ in Theorem 3.6 gives us the information about the cohomology ring of  the multiple weight variety of special type.

\begin{theorem}
Suppose that $G=SU(l+1)$ and $\lambda_{1}, \ldots, \lambda_{n} \in \mathfrak{t}^{*}_{++}$. Let $\mu \in \mathfrak{t}^{*}$ be a regular value of the moment map sufficiently close to $\lambda_{1}+\cdots+\lambda_{n}$, such that $\lambda_{1}+\cdots +\lambda_{n}-\mu$ is in the nice chamber. Then there exist cohomology classes $z_{1}, \ldots, z_{l} \in H^{2}(\mathcal{M}_{T};\mathbb{R})$  such that the cohomology ring $H^{*}(\mathcal{M}_{T};\mathbb{R})$ is given by
\[
H^{*}(\mathcal{M}_{T};\mathbb{R})\cong \frac{\mathbb{R}[z_{1},z_{2}, \ldots, z_{l}]}{(z_{l}^{n}, z_{l-1}^{n}(z_{l-1}+z_{l})^{n},\ldots, z_{1}^{n}(z_{1}+z_{2})^{n}\cdots (z_{1}+\cdots +z_{l})^{n})} .
\]
\end{theorem}

\begin{proof}
We set $\lambda_{i}=\sum\limits_{j=1}^{l}p_{i,j}\alpha_{j} \in \mathfrak{t}^{*}_{++}$ and $\mu=\sum\limits_{j=1}^{l}x_{j}\alpha_{j} \in \mathfrak{t}^{*}$ $(1\leq i \leq n,\ p_{i,j},x_{i} \in \mathbb{R})$. Let $\omega_{\mu}$ be the symplectic form on $\mathcal{M}_{T}$.

By Theorem 3.1, there exist cohomology classes $z_{i,j}$ and $v_{j} \in H^{2}(\mathcal{M}_{T};\mathbb{R})$ that generate $H^{*}(\mathcal{M}_{T};\mathbb{R})$ multiplicatively such that 
\[
[\omega_{\mu}]=\sum_{i=1}^{n}\sum_{j=1}^{l}p_{i,j}z_{i,j}+\sum_{j=1}^{l}x_{j}v_{j}.
\]
By Theorem 3.3, we have ${\rm vol}(\mathcal{M}_{T})=v_{l,n}(\lambda-\mu)$. Thus, ${\rm vol}(\mathcal{M}_{T})$ is a polynomial of $q_{1}=p_{1,1}+\cdots+p_{n,1}-x_{1},\ q_{2}=p_{1,2}+\cdots+p_{n,2}-x_{2},\ldots,\ q_{l}=p_{1,l}+\cdots +p_{n,l}-x_{l}$. By Lemma 2.6, we have the relations
\[
\left\{
\begin{array}{l}
z_{1,1}=z_{2,1}=\cdots=z_{n,1}=-v_{1}\\
z_{1,2}=z_{2,2}=\cdots=z_{n,2}=-v_{2}\\
\vdots\\
z_{1,l}=z_{2,l}=\cdots=z_{n,l}=-v_{l}.
\end{array}
\right.
\]
By Theorem 3.6, we have the system of differential equations
\[
\left\{
\begin{array}{l}
\partial^{n}_{l}v=0\\
\partial^{n}_{l-1}(\partial_{l-1}+\partial_{l})^{n}v=0\\
\vdots\\
\partial_{1}^{n}(\partial_{1}+\partial_{2})^{n}\cdots (\partial_{1}+\cdots +\partial_{l})^{n}v=0,
\end{array}
\right.
\]
where $\partial_{i}=\frac{\partial}{\partial q_{i}}$ for $1\leq i \leq l$.
By Lemma 2.5, we have the relations
\[
\left\{
\begin{array}{l}
z_{l}^{n}=0\\
z_{l-1}^{n}(z_{l-1}+z_{l})^{n}=0\\
\vdots\\
z_{1}^{n}(z_{1}+z_{2})^{n}\cdots (z_{1}+\cdots +z_{l})^{n}=0.
\end{array}
\right.
\]
By Theorem 3.8, these are all relations in $H^{*}(\mathcal{M}_{T};\mathbb{R})$. Thus, we have
\[
H^{*}(\mathcal{M}_{T};\mathbb{R})\cong \frac{\mathbb{R}[z_{1},z_{2}, \ldots, z_{l}]}{(z_{l}^{n}, z_{l-1}^{n}(z_{l-1}+z_{l})^{n},\ldots, z_{1}^{n}(z_{1}+z_{2})^{n}\cdots (z_{1}+\cdots +z_{l})^{n})} .
\]
\end{proof}

\begin{remark}
In Theorem 3.9, we consider the special case, where $\lambda-\mu$ is in the nice chamber. In [6], a presentation of the cohomology ring of $\mathcal{M}_{T}$ is given  for general $\mu$ in the case $G=SU(l+1)$ and $n=1$. It is interesting to compare the result in Theorem 3.9 with the presentation in [6].
\end{remark}

\begin{example}
We consider the case in Example 3.4 and $h:=\lambda_{1}+\lambda_{2}-\mu \in \mathfrak{c}_{\rm nice}= \mathbb{R}_{>0}(\alpha_{1}+\alpha_{2})+\mathbb{R}_{>0}\alpha_{2}$. Then, the symplectic ${\rm vol}(\mathcal{M}_{T})$ is given by
\[
{\rm vol}(\mathcal{M}_{T})=\frac{1}{12}(p+r-x)^{3}(-(p+r-x)+2(q+s-y)). 
\]
Note that ${\rm vol}(\mathcal{M}_{T})$ is a polynomial of $X:=p+r-x$ and $Y:=q+s-y$, and so we have 
\[
{\rm vol}(\mathcal{M}_{T})=\frac{1}{12}X^{3}(2Y-X).
\]
and satisfies the assumption in Theorem 3.3 and Theorem 3.6.
Therefore ${\rm vol}(\mathcal{M}_{T})$ satisfies the following differential equations
\[
\left\{
\begin{array}{l}
\partial^{2}_{Y}{\rm vol}(\mathcal{M}_{T})=0\\
\partial^{2}_{X}(\partial_{X}+\partial_{Y})^{2}{\rm vol}(\mathcal{M}_{T})=0.\\
\end{array}
\right.
\]

By Theorem 3.9, we have
\[
H^{*}(\mathcal{M}_{T};\mathbb{R})\cong \frac{\mathbb{R}[z_{1},z_{2}]}{(z_{2}^{2},z_{1}^{2}(z_{1}+z_{2})^{2})}.
\]
This presentation of the cohomology ring of $\mathcal{M}_{T}$ also gives ${\rm dim}\ H^{2}(\mathcal{M}_{T};\mathbb{R})$ $=2$, ${\rm dim}\ H^{4}(\mathcal{M}_{T};\mathbb{R})=2$, ${\rm dim}\ H^{6}(\mathcal{M}_{T};\mathbb{R})=2$ and ${\rm dim}\ H^{8}(\mathcal{M}_{T};\mathbb{R})=1$. 
\end{example}



\end{document}